\documentclass[11pt]{amsart}
\usepackage{amsmath,amssymb,amsthm}
 
\newcommand{\sortnoop}[1]{}

\theoremstyle{plain}
\newtheorem{theorem}{Theorem}
\newtheorem{proposition}[theorem]{Proposition}
\newtheorem{lemma}[theorem]{Lemma}
\newtheorem{corollary}[theorem]{Corollary}

\theoremstyle{definition}

\newtheorem*{theorem*}{Theorem}

\newcommand{\reff}[1]{(\ref{#1})}

\newcommand{\findemo}{\hfill $\square$\par
\vspace{0.3cm} 
\noindent}

\newcommand{\rw}{{\rm w}}

\newcommand{\rP}{{\rm P}}
\newcommand{\rE}{{\rm E}}
\newcommand{\PP}{{\rm \bf P}}

\newcommand{\ca}{{\mathcal A}}
\newcommand{\cb}{{\mathcal B}}

\newcommand{\ccc}{{\mathcal C}}

\newcommand{\ce}{{\mathcal E}}

\newcommand{\cg}{{\mathcal G}}
\newcommand{\ch}{{\mathcal H}}

\newcommand{\cw}{{\mathcal W}}

\newcommand{\E}{{\mathbb E}}

\newcommand{\N}{{\mathbb N}}
\renewcommand{\P}{{\mathbb P}}

\newcommand{\R}{{\mathbb R}}

\newcommand{\ind}{{\bf 1}}
\newcommand{\range}{{\mathcal R}}

\newcommand{\supp}{{\rm supp}\;}

\newcommand{\norm}[1]{\mathop{\parallel\! #1 \! \parallel}\nolimits}
\newcommand{\val}[1]{\mathop{\left| #1 \right|}\nolimits}
\newcommand{\inv}[1]{\mathop{\frac{1}{ #1}}\nolimits}
\newcommand{\expp}[1]{\mathop {\mathrm{e}^{ #1}}}
\newcommand{\capa}{\mathop{\mathrm{cap}}\nolimits}

\newcommand{\Capa}{\mathop{\mathrm{Cap}}\nolimits}
\newcommand{\capaD}{\mathop{\mathrm{cap_D}}\nolimits}
\newcommand{\normt}[1]{\mathop{\parallel\! #1 \! \parallel}\nolimits}

\begin{document}

\title[Characterization of G-regularity for super-Brownian motion]{Characterization of G-regularity for super-Brownian motion and consequences for parabolic partial differential equations}
\date{\today}
\author{Jean-Fran\c{c}ois Delmas}
\address{MSRI, 1000 Centennial Drive, Berkeley, CA 94720, U.S.A. \\
and ENPC-CERMICS, 6 av. Blaise Pascal, Champs-sur-Marne, 77455 Marne La Vall\'ee, France.}
\email{delmas@msri.org}
\thanks{The research of the first author was done at the \'Ecole Nationale des Ponts et Chauss\'ees and at MSRI, supported by NSF grant DMS-9701755.}
\author{Jean-St\'ephane Dhersin}
\address{UFR de Math\'ematiques et d'Informatique, Universit\'e Ren\'e Descartes, Paris, France}
\email{dhersin@math-info.univ-paris5.fr}

\begin{abstract}
We give a characterization of G-regularity for super-Brownian motion and the
Brownian snake. More precisely, we define a capacity on $E=(0,\infty )\times \R^d$, 
which is not invariant by translation. We then prove that the
hitting probability of a Borel set $A\subset E$ for the graph of the Brownian snake starting
at $(0,0)$ is comparable, up to multiplicative constants, to its
capacity. This implies that super-Brownian motion started at time $0$ at the Dirac mass 
$\delta_0$ hits immediately $A$ (that is $(0,0)$ is G-regular for $A^c$) 
if and only if its capacity is
infinite. As a direct consequence, if $Q\subset E$ is a domain such that $(0,0)\in \partial Q$, we give a necessary and sufficient condition for the existence  on $Q$  of a positive solution of $\partial_t u+\frac{1}{2}\Delta u =2u^2$ which blows up at $(0,0)$. We also give an estimation of the hitting probabilities for the support of super-Brownian motion at fixed time. We prove that if $d\geq 2$, the support of super-Brownian motion is intersection-equivalent to the range of Brownian motion.
\end{abstract}



\maketitle 

\section{Introduction}
The purpose of this paper is to give a characterization of the so called
G-regularity for super-Brownian motion introduced by Dynkin~\cite{d:spde}. Let us recall that a
point $(r,x)\in {\mathbb R}\times{\mathbb R}^d$ is G-regular for an open set
$Q\subset  {\mathbb R}\times{\mathbb R}^d$ if a.s. the graph of a super-Brownian motion started at time $r$ with the Dirac mass at $x$ immediately intersects $Q^c$, the complementary    of $Q$. (This definition can be extended to any Borel set.) We also recall that this is equivalent to the existence of nonnegative solutions of the equation $\frac{\partial u}{\partial t}+\frac{1}{2}\Delta u=2\,u^2$ on the open set $Q$, which blow up at $(r,x)\in \partial Q$.

Let $E=(0,\infty )\times{\mathbb R}^d$.
We prove that $(0,0)$ is
$G$-regular for a Borel set $A\subset \R\times \R^d$ if and only if the  capacity of $A^c\cap E$ is infinite, for the following capacity: for any Borel set $A'\subset E$, 
\[
\capa (A')=\left[\inf I(\nu)\right]^{-1}, \quad \text{where }
\]
\[
I(\nu)= \iint_E dsdy\; p(s,y)\left(\iint_E \nu(dt,dx)\frac{p(t-s,x-y) }{p(t,x)}\right)^2,
\]
and $p$ denotes the heat kernel:
\[
p(t,x)=\left\{
\begin{array}{ll}
   (2 \pi t)^{-d/2}\expp{-\val{x}^2/2t} \quad &\text{ if } (t,x)\in E,\\
   0 &\text{ if } (t,x)\in (-\infty ,0]\times \R^d.
\end{array}
\right.
\] 
($\val{\cdot}$ denotes the Euclidean norm on $\R^d$.)
The infimum is taken over all probability measures $\nu$ on $E$ such that $\nu(A')=1$. 
Notice this capacity  is not invariant by translation in time or space. This capacity arises naturally when one consider the Brownian snake, a useful tool to study super-Brownian motion. Indeed, using potential theory of symmetric Markov process, $I(\nu)$ can be viewed as the energy, with respect to the Brownian snake, of a certain probability measure (see section 4 for more details).

We extend a result due to Dhersin and Le~Gall~\cite{dlg:ktsbm} where the authors study $G$-regularity of $(0,0)$ for sets $Q=\{(s,y)\in E;\,|y|<\sqrt{s}\,h(s)\}$, where $h$ is a positive  decreasing function defined on $(0,\infty )$. Our result can also be viewed as a parabolic extension of the Wiener's test proved by Dhersin and Le~Gall~\cite{dlg:wtsbmbs} in an elliptic setting.

The proof of our results relies on the   Brownian snake introduced by Le~Gall. We only give definition and some properties for completeness in this paper, and refer to Le~Gall~\cite{lg:pvmp,lg:pvmppde} for a detailed presentation. We will use time inhomogeneous notations.

Let $(r,x)\in {\mathbb R}\times {\mathbb R}^d$ be a fixed point. We denote by ${\mathcal W}_{r,x}$ the set of all stopped paths in ${\mathbb R}^d$ started at $x$ at time $r$. An element ${\rm w}$ of ${\mathcal W}_{r,x}$ is a continuous mapping ${\rm w}:[r,\zeta]\to{\mathbb R}^d$ such that ${\rm w}(r)=x$, and $\zeta=\zeta_{(\rm w)}\in [r,\infty )$ is called its lifetime. We denote by $\hat \rw $ the end point $\rw(\zeta)$.  With the metric 
$d(\rw,\rw')=\val{\zeta_{(\rw)} -\zeta_{(\rw')}}+\sup_{s\geq
  r}\val{\rw(s\wedge\zeta_{(\rw)})-\rw'(s\wedge\zeta_{(\rw')})}$,
the space $\cw_{r,x}$ is a Polish space. The Brownian snake started at $x$ at time $r$ is a continuous strong Markov process $W=(W_s,s\geq 0)$ with values in $\cw_{r,x}$, whose law is characterized by the following two properties. 
\begin{itemize}
   \item[(i)] The lifetime process  $\zeta=\left(\zeta_s=\zeta_{(W_s)},s\geq 0\right)$ is a reflecting Brownian motion in $[r,\infty )$.
   \item[(ii)] Conditionally given $\left(\zeta_s,s\geq 0\right)$, the process $\left(W_s,s\geq 0\right)$   is a time-inhomogeneous continuous Markov process, such that for  $s'\geq s$:
\begin{itemize}
\item[$\cdot$] $W_{s'}(t)=W_s(t)$ for  $r\leq t\leq m(s,s')={{\displaystyle }
\inf_{v\in [s,s']}\zeta_v}$.
\item[$\cdot$] $(W_{s'}(m(s,s')+t)-W_{s'}(m(s,s')),0\leq t\leq
\zeta_{s'}-m(s,s'))$
is a Brownian motion in ${\mathbb R}^d$ independent of $W_s$.
\end{itemize}
\end{itemize}

From now on we shall consider the canonical realization of the process $W$ defined on the space $\Omega=C(\R^+,\cw_{r,x})$, and denote by $\ce_{\rw}$ the law of $W$ started at $\rw\in\cw_{r,x}$.  The trivial path $\textbf{x}_r$ such that $\zeta_{(\textbf{x}_r)}=r$, $\textbf{x}_r (r)=x$ is clearly a regular point for the process $\left(W, \ce_{\rw}\right)$. We denote by $\N_{r,x}$ the excursion measure outside $\left\{\textbf{x}_r\right\}$. Notice that $\N_{r,x}$ is an infinite measure. 
The distribution of $W$ under $\N_{r,x}$ can be characterized as above, except that now  the lifetime process $\zeta$ is distributed according to the It\^o measure of excursions of linear reflecting Brownian motion in $[r,\infty )$. We  normalize $\N_{r,x}$ so that, for every $\varepsilon>0$,
\[
\N_{r,x}\left[\sup_{s\geq 0}\zeta_s>\varepsilon+r\right]=\inv{2\varepsilon}.
\]
Let $\sigma=\inf\left\{s>0;\zeta_s=r\right\}$ denote  the duration of the
excursion of $\zeta$ under $\N_{r,x}$. The graph $\cg^*$ of $W$ is defined under
$\N_{r,x}$ by
\[
\cg^*=\left\{(t,W_s(t)); r< t\leq  \zeta_s, 0<  s <  \sigma\right\} =
\left\{(\zeta_s,\hat W_s); 0<s< \sigma \right\}.
\]

\noindent
We write $\cg^*(W)$ for $\cg^*$  when there is a risk of confusion.

Let us now explain the connection between the Brownian snake and super-Brownian motion. First of all, we introduce some notations. We denote by $(M_f,{\mathcal M}_f)$ the space of all finite measures on $\R^d$, endowed with the  topology of  weak  convergence.  We denote by  $\cb (S)$ (resp. $\cb_{b+}(S)$) the set of all real measurable  (resp. bounded nonnegative  measurable) functions defined on a polish space $S$. We also  denote by $\cb (S)$  the Borel $\sigma$-field  on $S$. For every measure $\nu\in M_f$,
and 
$f\in \cb_{b+}(\R^d)$, we shall write $(\nu,f)=\int f(y)\nu(dy)$. We also denote by $\supp \nu$ the closed support of the measure $\nu$. 

We consider under $\N_{r,x}$ the continuous version $\left(l^t_s,t> r,s\geq 0\right)$ of the local time of $\zeta$ at level $t$ and time $s$, and  define the  measure valued process $Y$ on $\R^d$ by setting  for every $t>r$, for every $\varphi\in \cb_{b+}(\R^d)$, 
\[
(Y_t,\varphi)=\int_0^\sigma dl_s^{t} \; \varphi(\hat W_s).
\]
Let $\cw_{r}=\bigcup _{x\in \R^d}\cw_{r,x}$.
Let $\mu$ be a finite measure on $\R^d$, and $\sum_{i\in I} \delta_{W^i}$ be a Poisson measure on
   $C(\R^+,\cw_{r})$ with intensity $\int\mu(dx)\N_{r,x}[\cdot]$. Then the process
   $X$ defined by $X_r=\mu$ and  $X_t=\sum_{i\in I}Y_t(W^i)$ if $ t>r$,
   is a super-Brownian motion started at time $r$ at $\mu$ (see  \cite{lg:pvmp,lg:pvmppde}).
We shall denote by $\PP_{r,\mu}$ (resp. $\PP_{r,x}$) the law of the super-Brownian motion started at time $r$ at $\mu$ (resp. at the Dirac mass $\delta_x$). We deduce from the normalization of $\N_{r,x}$ that, for every $t>r$, $\N_{r,x} \left[Y_t\neq 0\right]=1/2(t-r)<\infty $. This implies that  there is only a finite number of indices $i\in I$ such that $\cg^*(W^i)\cap [t,\infty )\times \R^d$ is non empty for $t>r$.

We consider the graph of $X$: \label{com:cgX}
\[
\cg(X)=\bigcup_{\varepsilon>r}\left(\overline{\bigcup_{t\geq \varepsilon} \{t\}\times \supp X_t }\right)=\bigcup _{i\in I} \cg^*(W^i),
\]
where $\bar A$ denotes the closure of $A$. A set $A\subset \cb(\R\times \R^d)$ is called $G$-polar if $\PP_{r,x} [\cg(X)\cap A\neq \emptyset]=0$ for every $(r,x)\in \R\times \R^d$. From Poisson measure theory, we have
\begin{equation*}
\PP_{r,x} [\cg(X)\cap A\neq \emptyset]=1-\expp{-\N_{r,x}[\cg^*\cap A\neq \emptyset]}.
\end{equation*}
Hence $A$ is $G$-polar if and only if $\N_{r,x}[\cg^*\cap A\neq \emptyset]=0$ for all $(r,x)\in \R\times \R^d$.  
We consider the capacity defined by: for $A\in \cb(\R\times \R^d)$,
\[
\capa' (A)=\left[\inf \iint_{\R\times \R^d} dsdy \left(\iint_{(s,\infty )\times \R^d}  \nu(dt,dx)p(t-s,x-y)\expp{-(t-s)/2 }\right)^2\right]^{-1},
\]
where the infimum is taken over all  probability measures $\nu$ on $\R\times \R^d$ such that $\nu(A)=1$. Dynkin proved (see Theorem 3.2 in \cite{d:sap}) that $A\in \cb(\R\times \R^d)$ is $G$-polar if and only if $\capa' (A)=0$. (We have $\capa' (A)=0 \Leftrightarrow \N_{r,x}[\cg^*\cap A\neq \emptyset]=0$ for all $(r,x)\in \R\times \R^d$.) It is easy to check that if $A\subset E$ is a compact set then
\[
\capa' (A)=0 \Leftrightarrow \capa (A)=0.
\]
This can be extended to all Borel subsets of $E$ since the two capacities are inner capacities (see Meyers \cite{m:tcpf}).  
In fact it seems  more relevant to consider the capacity $\capa$ to characterize $G$-regularity, as we shall see. We have the following quantitative theorem.
\begin{theorem}
   \label{th:ucap}
There exists a constant $C_0$ such that for any $A\in \cb(E)$,
\[
4^{-1} \capa(A)\leq \N_{0,0}[\cg^*\cap A\neq \emptyset]\leq  C_0 \capa(A).
\]
\end{theorem}

The proof of Theorem~\ref{th:ucap} is split in two parts. In section 2, we introduce a capacity associated with a weighted Sobolev space, which is equivalent to the capacity $\capa$. In section 3, using the connections between super-Brownian motion and partial differential equations, we prove the upper bound with this new capacity, and hence for the capacity $\capa$. The lower bound is obtain in section 4, by using additive functionals of the Brownian snake introduced in \cite{dlg:wtsbmbs}.

Now, for  $A\in \cb( \R\times \R^d)$, we consider under $\PP_{r,x}$ the random time
\[
\tau_A=\inf\{t>r, \left(\{t\}\times \supp X_t \right)\cap A\neq\emptyset\}.
\]
Arguments similar as those of  \cite{dlg:wtsbmbs} yield that  $\tau_A$ is a stopping time for the natural filtration of $X$ completed the usual way. Thus we have $\PP_{r,x}(\tau_A=r)=1$ or $0$.
\label{com:ta}
Following Dynkin~\cite[section II-6]{d:spde}, we say a point $(r,x)\in  \R\times \R^d$ is $G$-regular for $A^{c}$ if $\PP_{r,x}$-a.s. $\tau_A=r$. Let $A^{Gr}$ denote the set of all points that are G-regular for $A^{c}$.  From the known path properties of super-Brownian motion it is obvious that int$(A)\subset A^{Gr}\subset \bar A$, where int$(A)$ denotes the interior of $A$. We set $T_A=\inf \left\{s>0, (\zeta_s,\hat W_s) \in A\right\}$. Following \cite{dlg:wtsbmbs} it is easy to deduce from Theorem~\ref{th:ucap} the next result.
\begin{proposition}
\label{prop:sreg}
   Let $A\in \cb(\R\times \R^d)$. The following properties are equivalent:
\begin{enumerate}
   \item $(r,x)$ is $G$-regular for $A^{c}$;
   \item $\N_{r,x}\left[\cg^*\cap A\neq\emptyset \right]=\infty $;
   \item $\ce_{\mbox{\rm{$\textbf{x}_r$}}}$-a.s. $T_A=0$;
   \item $\capa (A_{r,x}\cap E)=\infty $, where $A_{r,x}=\{(s,y); (s+r,y+x)\in A\}$.
\end{enumerate}
\end{proposition}

We can give a straightforward analytic consequence of
Proposition~\ref{prop:sreg} and the link between super-Brownian motion
and nonlinear differential equation. 
\begin{corollary}
   Let $Q$ a domain in $E$ such that $(0,0)\in \partial Q$. The following three conditions are equivalent.
\begin{enumerate}
   \item $(0,0)$ is $G$-regular for $Q$ ;
   \item $\capa (Q^c\cap E)=\infty$ ;
   \item There exists a nonnegative solution of $\frac{\partial u}{\partial t} +\frac{\Delta}{2} u=2 u^2 $ in $Q$ such that 
\[
\lim_{(s,y)\rightarrow (0,0),\; (s,y)\in Q} u(s,y)=\infty.
\]
\end{enumerate}
\end{corollary}

The equivalence of assertions 1) and 3) is due to Dynkin~\cite[Theorem II.6.1]{d:spde}.  The equivalence of 1) and 2) is given by Proposition~\ref{prop:sreg}.

Finally, using Theorem~\ref{th:ucap} we give  in section 5 an estimation of the hitting probability of the support of $X_1$. And we prove that in dimension $d\geq 2$, the support of super-Brownian motion and the range of $d$-dimensional Brownian motion are intersection-equivalent.

\section{Equivalence of capacities for a weighted Sobolev space}
\label{sec:sobolev}

\noindent
In this section, we introduce a new capacity, associated with a weighted Sobolev space, which is equivalent to the capacity $\capa$. This capacity will be very useful in the next section to  prove the upper bound for Theorem~\ref{th:ucap}.

If $S$ is an open subset of $\R^p$, we denote by  $C_0^\infty (S)$  the set of all  functions of class $C^\infty $ defined on $S$ with compact support. If $f$ is a measurable function defined on $S$ then $\norm{f}_\infty =\sup_{x\in S} \val{f(s)}$. 
We consider the Hilbert space $L^2(p)=\left\{f \in \cb(E);\; \norm{f}_{(p)}<\infty  \right\}$,
where $\norm{f}_{(p)}^2=\iint_E dtdx\; p(t,x) f(t,x)^2$. 

Notice the kernel defined on $E\times E$ by $k(t,x;s,y)=p(t-s,x-y)p(t,x)^{-1}$
is nonnegative  and lower semi-continuous. 
Thus we can introduce the operator $\Lambda$ defined on the set of nonnegative functions  $f\in \mathcal{B}(E)$ by:
\[
\Lambda(f)=p^{-1} [p*(pf)]=\iint_E dsdy\; k(\cdot,\cdot;s,y)p(s,y) f(s,y),
\]
where $*$ denotes the usual convolution product on $E$. Furthermore, the function $\Lambda(f)$ is even lower semi-continuous (see  \cite[Lemma 2.2.1]{f:tplcs}).

We define the  capacity $\Capa$ on $E$ in the following way: if $A\subset E$, then
\[
\Capa (A)=\inf \left\{ \norm{f}_{(p)}^2; \;f\geq 0,\; f\in L^2(p),\; \Lambda(f)\geq 1 \quad\text{on} \quad A \right\},
\]
with the convention $\inf \emptyset=\infty $.
Notice this capacity is not invariant by translation in time or space.
This capacity is an outer capacity (see Meyers \cite[Theorem 1]{m:tcpf}).
Moreover, it coincides with the capacity  $\capa$ on the analytic sets (see \cite[Theorem 14]{m:tcpf}). Now, we want to connect this capacity to an analytic capacity (see Baras and Pierre \cite{bp:ppsdm} for similar results but with different norms). Therefore we  consider the weighted Sobolev space $W_D$ which is the completion of $C^\infty_0(E)$ with respect to the norm $\norm{\cdot}_D$, defined by 
\[
\norm{\varphi}_D^2=\norm{\partial_t \varphi}_{(p)}^2+ \sum_{i=1}^d \norm{\partial_i (\log p) \; \partial_i \varphi }_{(p)}^2+ \sum_{i=1}^d \norm{\partial^2_{ii}\varphi}_{(p)}^2,\quad \varphi\in {C}_{0}^{\infty }(E),
\]
with the usual notations $\partial_t g(t,x)=\frac{\partial g}{\partial t} (t,x)$,  $\partial_i g(t, x)=\frac{\partial g}{\partial x_i} (t,x)$ for $x=(x_1,\cdots,x_d)\in \R^d$ and $\partial^2_{ii}=\partial_i\partial_i$.
Notice the non zero constants do not belong to $W_D$. We can introduce the  outer capacity $\capaD$ associated to $W_D$  defined as follows.  For any compact set $K\subset E$,  we set 
\begin{align*}
   \capaD (K)
&=\inf \left\{\norm{\varphi}_D^2; \;\varphi\in C_0^\infty (E), \;\varphi\geq 0, \;\varphi\geq 1 \text{ on } K \right\}\\
&=\inf \left\{\norm{\varphi}_D^2; \;\varphi\in C_0^\infty (E), \;\varphi\geq 0, \;\varphi\geq 1 \text{ on a neighborhood of } K \right\}.
\end{align*}
Then we set for any open set $G\subset E$, 
\begin{equation}
   \label{eq:capdGK}
\capaD (G)=\sup\left\{\capaD(K); \;K\subset G, \;K \quad\text{compact} \right\},
\end{equation}
and, for  any analytic set $A\subset E$,
\[
\capaD (A)=\inf\left\{\capaD(G); A\subset G, \;G \quad\text{open} \right\}.
\]
Notice the definition is consistent (see \cite{ah:fspt} for example).
\begin{proposition}
There exists a constant $C$ such that for any set $A\subset E$,
\[
\Capa (A)\leq \capaD(A)\leq C \Capa(A).
\]
\end{proposition}

\noindent
\textbf{Proof.}
Since the two capacities are outer capacity, it is enough to consider open sets. Now, using  \reff{eq:capdGK} and \cite[Theorem 8]{m:tcpf}, we see it is enough to consider compact sets.

Let us introduce the operator $H={\partial}_t-\inv{2}{\Delta}$. We consider a non empty compact set $K\subset E$. Let  $\varphi\in C_0^\infty(E) $ be such that $\varphi\geq 0$ and $\varphi\geq 1$ on $K$. Notice that (in the distribution sense) $ H p=\delta_{(0,0)}$, where $\delta_{(0,0)}$ is the Dirac mass at ${(0,0)}\in {\mathbb R}\times\R^{d}$. Then we have $p* [H (p\varphi)]= (H p)*(p\varphi)= p\varphi$. The function $f=p^{-1} \val{H(p\varphi)}=\val{H(\varphi)-(\nabla \log p, \nabla \varphi)}$ is nonnegative and $\Lambda(f)\geq \varphi$. Thus we have $\Lambda(f)\geq 1$ on $K$. We also have  $\norm{f}_{(p)}\leq \norm{\varphi}_D$. Hence we have $\Capa (K)\leq \norm{\varphi}_D^2$. The first inequality follows by taking a sequence $(\varphi_n)$ such that $\norm {\varphi_n}_D^2$ converges to $\capaD (K)$.

To prove the other inequality, let us consider a nonnegative function $f_1\in L^2(p)$,  such that $\Lambda (f_1)\geq 1$ on $K$. Notice it implies $\norm{f_1}_{(p)}>0$. Let $\delta >0$. It is easy to construct a function $\varepsilon\in L^2(p)$ such that $\varepsilon >0$ on $E$ and $\norm{\varepsilon}_{(p)}\leq \delta\norm{f_1}_{(p)}$. We set $f_2=f_1+\varepsilon$. Since the function $\Lambda (f_2)$ is lower semi-continuous, the set $\{(t,x)\in E; \;\Lambda (f_2)>1\}$  is open and it also contains $K$. It is then obvious that for $\delta'>0$  small enough, if we set $f_3(t,x)=f_2(t,x)\ind_{\{\delta'< t<{\delta'}^{-1}, \val{x}<{\delta'}^{-1} \}}$ for $(t,x)\in E$, we get $\Lambda (f_3)> 1$ on an open set containing $K$. 
Let us introduce a nonnegative function $h\in C_0^\infty (E)$  such that $\iint_E h(t,x)dtdx=1$. For $\theta>0$, we write $h_\theta(t,x)=\theta^{-d-1} h(t/\theta,x/\theta)$. Now using the uniform continuity of $p$ on $[\delta'/2,\infty )\times \R^d$, it is easy to see that if $f=h_\theta *f_3$ , then $\Lambda (f)> 1$ on an open set containing $K$ for $\theta$ small enough. The function $f$ is nonnegative, belongs to $C_0^\infty(E) $ and the function $\Lambda(f)$ is of class $C^\infty $. We can choose $\delta$ and $\theta$ small enough so that $\norm{f}_{(p)}\leq 2\norm{f_1}_{(p)}$.
\label{com:convol}

Let $\alpha \in C_0^\infty ([0,\infty ))$ such that $0\leq \alpha\leq 1$, $\alpha=1$ on $[0,1/2]$ and $\alpha=0$ on $[1,\infty) $. Let $\xi \in C_0^\infty (\R^d)$ such that $0\leq \xi\leq 1$ and $\xi=1$ in a neighborhood of $0$. We define $\alpha_n(t)=\alpha (t/n)$ and $\xi_n(x)=\xi(x/n)$. The function $\varphi_{n}= \alpha_n \xi_n \Lambda(f)$ belongs to  $C_0^\infty (E)$, is nonnegative and $\varphi_{n} \geq 1$ on a neighborhood of $K$ for $n$ great enough. 

Let us now give two key lemmata.  If $M$ is a bounded operator from $L^2(p)$ into itself, we denotes by  $\normt{M}_{(p)}= \sup\{\norm{M(f)}_{(p)}; f\in L^2(p),\; \norm{f}_{(p)}=1\}$ its norm. We define the operator $\Lambda_0$: for $f\in \cb(E)$ nonnegative, $\Lambda_0(f)(t,x)=t^{-1} \Lambda(f)(t,x)$, $(t,x)\in E$. For $T>0$, let us introduce $E_T=(0, T)\times \R^d$. 

\begin{lemma}
\label{lem:opL}
   The operators $\ind_{E_T}\Lambda$ and $\Lambda_0$ are bounded operators from $L^2(p)$ into itself. Furthermore, we have 
$\normt{\ind_{E_T}\Lambda}_{(p)}\leq T/\sqrt{2}$ and $\normt{\Lambda_0}_{(p)}\leq {2}$.
\end{lemma}
\noindent
\textbf{Proof} of Lemma \ref{lem:opL}. 
Let $f\in L^2(p)$. The Cauchy-Schwarz inequality implies:
\begin{align*}
\norm{\Lambda_0(f) }_{(p)}^2&= \iint_{E} dt dx\; t^{-2} p(t,x)^{-1} \left[\iint_E  ds dy\; p(t-s,x-y) p(s,y) f(s,y) \right]^2\\
&\leq \iint_{E} dt dx\; t^{-2}p(t,x)^{-1} \iint_E   ds' dy'\; p(t-s',x-y') p(s',y'){s'}^{-1/2}\ind_{s'\leq t} \\
&\phantom{\leq \iint_{E_T} dt dx\; t^{-2}p(t,x)^{-1} }
\iint_E  ds dy\; p(t-s,x-y) p(s,y) {s}^{1/2}f(s,y)^2 \\
& =2\iint_{E} dt dx\; t^{-2}  t^{1/2}\iint_E  ds dy\; p(t-s,x-y) p(s,y) {s}^{1/2}f(s,y)^2 \\
& \leq 2 \iint_{E}  ds dy\; p(s,y) f(s,y)^2 {s}^{1/2}\int_s^\infty  t^{-3/2}dt \leq 4\norm{f}^2_{(p)}.
\end{align*}
Hence the operator  $\Lambda_0$ is a bounded operator from $L^2(p)$ into itself. And we have $\normt{\Lambda_0}_{(p)}\leq {2}$. The operator $\ind_{E_T}\Lambda$ can handled in a very similar way.
\findemo


\begin{lemma}
\label{lem:opL1}
   The operators defined on $C_0^\infty (E) $ by: $g \in C_0^\infty(E)  $
\begin{align*}
   \Lambda_1(g)& ={\partial_t} \Lambda(g),\\
\text{for}\quad i\in \{1,\cdots,d\}, \quad \Lambda_{2,i}(g)& =\inv{2}\partial^2_{ii} \Lambda(g),\\
\text{ and for}\quad i\in \{1,\cdots,d\}, \quad   \Lambda_{3,i}(g)& =\partial_{i}( \log p) \;\partial_{i}\Lambda(g),
\end{align*}
can be uniquely extended into bounded operators from $L^2(p)$ into itself. And we have
\begin{align}
\label{eq:opL1}   \normt{\Lambda_1}_{(p)}& \leq 1+3d,\\
\label{eq:opL2i}  \text{for}\quad i\in \{1,\cdots,d\}, \quad   \normt{\Lambda_{2,i}}_{(p)}& \leq 1,\\
\label{eq:opL3}  \text{for}\quad i\in \{1,\cdots,d\}, \quad   \normt{\Lambda_{3,i}}_{(p)}& \leq 4.
\end{align}
\end{lemma}

\noindent
The proof of this lemma is given in appendix.

We  now bound $\norm{\varphi_{n}}_D$.
Lemma \ref{lem:opL} provides an upper bound for $\norm{\partial_t \varphi_{n}}_{(p)}$:
\begin{align}
\nonumber
\norm{\partial_t \varphi_{n}}_{(p)}
& \leq \norm{\partial_t \alpha_n}_\infty \norm{\ind_{E_n}\Lambda(f)}_{(p)}+ \norm{\Lambda_1{f}}_{(p)}\\
\label{eq:dtphin}
& \leq \left(\norm{\partial_t \alpha}_\infty 2^{-1/2} +\normt{\Lambda_1}_{(p)}\right)\norm{f}_{(p)}.
\end{align}

\noindent
Using Lemma \ref{lem:opL} we derive an upper bound for $\sum_{i=1}^d  \norm{\partial_{i} \log p \;\partial_{i}\varphi_{n}}_{(p)}$:
\begin{align}
\nonumber
\sum_{i=1}^d  \norm{\partial_{i} \log p \;\partial_{i}\varphi_{n}}_{(p)}
&\leq  \sum_{i=1}^d\left(\norm{\Lambda_{3,i}(f)}_{(p)}+\sup_{x\in \R^d} \val{x_i \partial_i \xi(x)}\norm{\Lambda_0(f)}_{(p)}\right)\\
\label{eq:dilogpdiphin}
& \leq \sum_{i=1}^d \left(\normt{\Lambda_{3,i}}_{(p)}+\sup_{x\in \R^d} \val{x_i \partial_i \xi(x)}\normt{\Lambda_0}_{(p)} \right)\norm{f}_{(p)}.
\end{align}
In order to give an upper bound for $\sum_{i=1}^d \norm{\partial^2_{ii} \varphi_{  n}}_{(p)}$, we need an intermediary lemma.
\begin{lemma}
\label{lem:oPL2}
   There exists a constant $c_1$ (depending on $\xi$) such that for all $n\geq 1$, $g\in C_0^\infty (E)$, $i\in\{1,\cdots,d\}$,
\begin{equation*}
\norm{\ind_{E_n} \partial_i \xi_n \;\partial_i \Lambda(g) }_{(p)}\leq c_1 n^{-1/2}\norm{g}_{(p)}.
\end{equation*}
\end{lemma}

\noindent
\textbf{Proof.}
Recall that $\xi_n$ has compact support. Then, an integration by parts, Cauchy-Schwarz inequalities and Lemma~\ref{lem:opL}  give for $1\leq i\leq d$, 
\begin{align*}
   \norm{\ind_{E_n} \partial_i \xi_n \partial_i \Lambda(g)}_{(p)}^2
= &-\iint_E p \ind_{E_n} \Lambda(g)(\partial_i \xi_n )^2\partial^2_{ii} \Lambda(g)\\
&\hspace{1cm}- \iint_{E} p\ind_{E_n} \Lambda(g) (\partial_i \xi_n)^2 \partial_i \Lambda(g) \partial_i \log p\\
&\hspace{1cm}- 2\iint_{E} p\ind_{E_n} \Lambda(g) \partial_i \xi_n \partial_i \Lambda(g) \partial^2_{ii} \xi_n\\
\leq & \norm{\partial_i\xi_n}_\infty ^2\norm{\ind_{E_n} \Lambda(g)}_{(p)}\norm{\partial^2_{ii} \Lambda(g)}_{(p)}\\
&\hspace{1cm}+ \norm{\partial_i\xi_n}_\infty ^2\norm{\ind_{E_n} \Lambda(g)}_{(p)}\norm{\partial_i \Lambda(g)\;\partial_i \log p}_{(p)}\\
&\hspace{1cm}+2\norm{\partial^2_{ii} \xi_n}_\infty  \norm{\ind_{E_n} \Lambda(g)}_{(p)}\norm{\ind_{E_n} \partial_i \xi_n \partial_i \Lambda(g)}_{(p)}\\
\leq &2^{-1/2} n^{-1}\left[2\normt{\Lambda_{2,i}}_{(p)} +\normt{\Lambda_{3,i}}_{(p)}\right]\norm{\partial_i \xi}_\infty ^2\norm{g}_{(p)}^2\\
&\hspace{1cm}+2^{1/2} n^{-1}\norm{\partial^2_{ii} \xi}_\infty  \norm{g}_{(p)}\norm{\ind_{E_n} \partial_i \xi_n \partial_i \Lambda(g)}_{(p)}.
\end{align*}
Notice that if $a,b,c$ are positive then $a^2\leq c^2+ba$ implies $a\leq c+b$.
Thus we get
\begin{multline*}
\norm{\ind_{E_n} \partial_i \xi_n \partial_i \Lambda(g)}_{(p)}\\
\leq  2^{-1/4}n^{-1/2}\left[2\normt{\Lambda_{2,i}}_{(p)} +\normt{\Lambda_{3,i}}_{(p)}\right]^{1/2}\norm{\partial_i \xi}_\infty \norm{g}_{(p)}+  2^{1/2}  n^{-1}\norm{\partial _{i,i}^2  \xi}_\infty \norm{g}_{(p)},
\end{multline*}   
which ends the proof.
\findemo

Using this lemma and  Lemma \ref{lem:opL}, we get that
\begin{align}
\nonumber
  \sum_{i=1}^d \norm{\partial^2_{ii} \varphi_{  n}}_{p}
&\leq \sum_{i=1}^d\left[2\norm{\Lambda_{2,i}(f)}_{(p)}+\norm{\partial^2_{ii} \xi_n}_\infty \norm{\ind_{E_n} \Lambda(f)}_{(p)}+2\norm{\ind_{E_n}\partial_{i}\xi_n\; \partial_{i}\Lambda(f))}_{(p)}\right]\\
\label{eq:diiphin}
&\leq \sum_{i=1}^d \left[2\normt{\Lambda_{2,i}}_{(p)}+2^{-1/2}   n^{-1}\norm{\partial^2_{ii} \xi}_\infty +2c_1n^{-1/2}\right]\norm{f}_{(p)}.
\end{align}

\noindent
Then we deduce from \reff{eq:dtphin}, \reff{eq:dilogpdiphin}, \reff{eq:diiphin} and Lemma~\ref{lem:opL1} that there exists a constant $c_2$ independent of $f$ and $n\geq 1$ such that 
\[
\norm{\varphi_{n}}_D\leq c_2 \norm{f}_{(p)}.
\]
Thus we have $\norm{\varphi_n}_D\leq 2 c_2 \norm{f_1}_{(p)}$. The second inequality of the proposition is then obvious   with $C=4\,c_2^2$.
\findemo

We shall need the following lemma.
\begin{lemma}
\label{lem:phiK}
   There exists a constant $\gamma$ such that if $K\subset E$ is a  compact set  such that $\capaD (K)>0$, there exists  $\varphi\in C_0^\infty (E)$ which satisfies:

1) $0\leq \varphi\leq 1$,

2) $\varphi=1$ on a neighborhood of $K$,

3) $\norm{\varphi}_D^2\leq \gamma \capaD (K)$.
\end{lemma}
\noindent
The proof is classic, but we give it for completeness.

\noindent
\textbf{Proof.} 
Let $h\in C^\infty ([0,\infty ))$ such that $0\leq h\leq 1$, $h=0$ on $[0,1/4]$ and $h=1$ on $[3/4,\infty )$. Since $\capaD (K)>0$, there exists $g \in C_0^\infty (E)$ such that $g\geq 0$, $g\geq 1$ in a neighborhood of $K$, and $2\capaD (K)\geq \norm{g}_D^2$. Let $\varphi= h\circ g$. The function $\varphi\in C_0^\infty (E)$ satisfies 1) and 2). Let us check 3). We have
\begin{align*}
\norm{\partial_t \varphi}_{(p)}
&\leq \norm{h'}_\infty \norm{\partial_t g}_{(p)},\\
\norm{\partial_{i} \log p \;\partial_{i}\varphi}_{(p)}
&\leq \norm{h'}_\infty \norm{\partial_{i} \log p \;\partial_{i}g}_{(p)}\\
\norm{ \partial^2_{ii} \varphi}_{(p)}
&\leq \norm{h'}_\infty \norm{\partial^2_{ii} g}_{(p)}
+\norm{(h''\circ g) \;(\partial_i g)^2}_{(p)}.
\end{align*}
Only the upper bound for the  second right hand-side term of the last inequality is not obvious. We first search  an upper bound for $\norm{(\partial_i \varphi_1)^2/(1+\varphi_1)}_{(p)}$, where  $\varphi_1\in C^\infty _0(E)$ is a nonnegative function. An integration by parts  and Cauchy-Schwarz inequality give
\begin{align*}
   \iint_E p\frac{(\partial_i \varphi_1)^4}{(1+\varphi_1)^2}
&= 3 \iint_E p\; \partial_{ii}^2 \varphi_1 \frac{(\partial_i \varphi_1)^2}{1+\varphi_1}
+\iint_E p \;\partial_i \log p \;\partial_i \varphi_1 \frac{(\partial_i \varphi_1)^2}{1+\varphi_1}\\
&\leq (3 \norm{ \partial^2_{ii} \varphi_1}_{(p)} 
+\norm{\partial_{i} \log p \;\partial_{i}\varphi_1}_{(p)})\norm{{(\partial_i \varphi_1)^2}/{(1+\varphi_1)}}_{(p)}.
\end{align*}
Thus we get
\begin{equation}
\label{eq:maj-normediphi}
\norm{(\partial_i \varphi_1)^2/(1+\varphi_1)}_{(p)}\leq  3 \left(\norm{ \partial^2_{ii} \varphi_1}_{(p)}+ \norm{\partial_{i} \log p \;\partial_{i}\varphi_1}_{(p)}\right).
\end{equation}
Since we have $\val{h''(t)}\leq 2(1+t)^{-1} \norm{h''}_\infty $, taking $\varphi_1=g$ in the above inequality we deduce that
\begin{align*}
\norm{(h''\circ g) \;(\partial_i g)^2 }_{(p)}& \leq 2\norm{(\partial_i g)^2/(1+g)}_{(p)}\norm{h''}_{\infty }\\
& \leq 6(\norm{ \partial^2_{ii} g}_{(p)}+ \norm{\partial_{i} \log p \;\partial_{i}g}_{(p)})\norm{h''}_{\infty }.
\end{align*}
The previous inequalities imply there exists a constant $c$ depending only on $h$ and $d$ such that $\norm{\varphi}_D\leq c \norm{g}_D$. Thus 3) holds with $\gamma=2c^2$.
\findemo

\section{Upper bound for hitting probabilities} 
\label{sec:upper}

\noindent
In this section we prove the second inequality of Theorem \ref{th:ucap} for compact sets. Let us introduce $K\subset E_T$  a compact set such that $\capaD (K)>0$. Let $\varphi$ be as in Lemma \ref{lem:phiK}. We set $\varphi=0$ outside $E$. We introduce the function $\psi=1-\varphi$, which takes values in $[0,1]$. We consider the function $u$ defined on $\R\times \R^d$ by $u(t,x)=\N_{t,x}[\cg^*\cap K\neq\emptyset]$ ($\in [0,\infty ]$). With the convention $0.\infty =0$, the function $u\psi^4$ is bounded nonnegative and  of class $C^\infty $ on $\R\times \R^d$. Let $(B_t,t\geq 0)$ denote under $\rP_0$ a $d$-dimensional Brownian motion started from 0. It\^o's formula implies that for all $t\geq 0$, $\rP_0$-a.s., 
\begin{multline*}
   u\psi^4(t,B_t)=u\psi^4(0,0)+\int_0^t {\partial_t (u\psi^4)}(s,B_s) ds \\
+\int_0^t \frac{\Delta}{2} (u\psi^4)(s,B_s) ds +\int_0^t \nabla (u\psi^4)(s,B_s) dB_s.
\end{multline*}
Consider the stopping time $T_a=T\wedge \inf\{t>0; \val{B_t}\geq a\}$. We can then apply the stopping Theorem at time $T_a$ and get
\begin{align*}
\rE_0 u\psi^4(T_a,B_{T_a}) 
&=u(0,0)
+\rE_0 \int_0^{T_a} {\partial_t (u\psi^4)}(s,B_s) ds 
+\rE_0\int_0^{T_a} \frac{\Delta}{2} (u\psi^4)(s,B_s) ds\\
&=u(0,0)
+\rE_0 \int_0^{T_a} \bigg [2u^2 \psi^4+4u\psi^3{\partial_t \psi}+4 (\nabla u,\nabla \psi)\psi^2 \\
&\phantom{=u(0,0)+\E_0 \int_0^{T_a} [.]}
+6 u \psi^2 (\nabla \psi, \nabla \psi)+2 u \psi^3\Delta\psi  \bigg](s,B_s) ds.\\
\end{align*}
We have used that $\partial_t u+\inv{2} \Delta u =2u^2$ to get the last equality. Notice  that each integrand is either nonnegative  or bounded. 
By dominated convergence and monotone convergence, we get as $a$ goes to infinity
\begin{align*}
 u(0,0)+2\norm{u\psi^2\,\ind_{E_T}}_{(p)}^2=\rE_0 u\psi^4(T,B_{T}) 
- \iint_{E_T} p\bigg[
&4u\psi^3{\partial_t \psi}+4 (\nabla u,\nabla \psi)\psi^2 \\
&+6 u \psi^2 (\nabla \psi, \nabla \psi)+2 u \psi^3\Delta\psi  \bigg].
\end{align*}
Since $K\subset E_T$, we deduce that $u(t,x)=0$ for $t\geq T$. Thus we have:
\begin{align}
\nonumber
 u(0,0)+2\norm{u\psi^2}_{(p)}^2=
- \iint_{E} p\bigg[
&4u\psi^3{\partial_t \psi}+4 (\nabla u,\nabla \psi)\psi^2 \\
\label{eq:egalu00}
&+6 u \psi^2 (\nabla \psi, \nabla \psi)+2 u \psi^3\Delta\psi  \bigg].
\end{align}
We now  bound the right hand side. Using  Cauchy-Schwarz inequality, that $0\leq \psi\leq 1$, and that $-\varphi$ and $\psi$ have the same derivatives, we get
\begin{align*}
   -\iint_E pu\psi^3 {\partial_t \psi} &\leq \norm{u\psi^2}_{(p)} \norm{{\partial_t \varphi }}_{(p)},\\
 -\iint_E pu \psi^3 \partial^2_{ii} \psi &\leq \norm{u\psi^2}_{(p)} \norm{\partial^2_{ii}\varphi}_{(p)},
\end{align*}
\begin{align*}
\text{and}\quad -\iint_E p u\psi^2 (\nabla \psi, \nabla \psi)
& \leq   \norm{u\psi^2}_{(p)} \sum_{i=1}^d \norm{(\partial_i\varphi)^2}_{(p)}\\
& \leq   2\norm{u\psi^2}_{(p)} \sum_{i=1}^d \norm{{(\partial_i\varphi)^2}/{(1+\varphi)}}_{(p)}\\
& \leq   6\norm{u\psi^2}_{(p)} \sum_{i=1}^d (\norm{ \partial^2_{ii} \varphi}_{(p)}+ \norm{\partial_{i} \log p \;\partial_{i}\varphi}_{(p)}),\\
\end{align*}
where we have used  \reff{eq:maj-normediphi} with $\varphi_1=\varphi$ for the last inequality. 
Now an integration by parts and Cauchy-Schwarz inequality give
\begin{align*}
   -\iint_E p \psi^3 (\nabla u, \nabla \psi)
&= \iint_E pu \psi^2\left[\psi (\nabla \log p, \nabla \psi)+3 (\nabla \psi, \nabla \psi)+\psi \Delta\psi\right]\\
&\leq  \norm{u\psi^2}_{(p)} \sum_{i=1}^d \left[ \norm{\partial_{i} \log p \;\partial_{i}\varphi}_{(p)}
+3\norm{ (\partial_{i} \varphi)^2}_{(p)}+\norm{ \partial^2_{ii} \varphi}_{(p)} \right]\\
&\leq  19 \norm{u\psi^2}_{(p)} \sum_{i=1}^d \left[ \norm{\partial_{i} \log p \;\partial_{i}\varphi}_{(p)}
+\norm{ \partial^2_{ii} \varphi}_{(p)} \right],
\end{align*}
where we have used  again \reff{eq:maj-normediphi} for the last inequality. 
Taking those results together, we deduce from \reff{eq:egalu00} that 
\[
u(0,0)+2\norm{u\psi^2}_{(p)}^2\leq  c_3 \norm{u\psi^2}_{(p)}\norm{\varphi}_D,
\]
where the constant $c_3$ depends only on $d$. Since $\norm{u\psi^2}_{(p)}$ is finite (recall $u\psi$ is bounded, and zero on $[T,\infty )\times \R^d$), this implies that $\norm{u\psi^2}_{(p)}\leq  c_3 \norm{\varphi}_D$ and hence $u(0,0)\leq c_3^2\norm{\varphi}_D^2$. This last inequality and the definition of $\varphi$ imply that
\[
\N_{0,0}[\cg^*\cap K\neq\emptyset]=u(0,0)\leq c_3^2\gamma  \capaD(K)\leq  c_3^2 \gamma C \Capa (K) = c_3^2 \gamma C \capa (K).
\]

\section{Lower bound for hitting probabilities and proof of Theorem \ref{th:ucap}}
\label{sec:lower}

\noindent
In this section, we prove the first inequality  of Theorem \ref{th:ucap} for
compact sets. Let us introduce a compact set $K\subset E$, $\nu$ a
probability measure on $K$, and $T>0$ such that  $K\subset E_T$. We consider the probability measure $\mu$ defined on $\cw_{0,0}$ by
\[
\mu(d\rw) =\iint_E \nu(dt,dx) \; {\rm P}_0^{t,x}(d\rw),
\]
where ${\rP}_0^{t,x}$ is the law on $\cw_{0,0}$ of the Brownian bridge starting at time $0$ at point $0$ and ending at time $t$ at point $x$. Notice that the measure $\mu$ is in fact a measure on $\cw_{0,0}^*$, the set of non trivial path in $\cw_{0,0}$ (a trivial path is a path of lifetime zero). The measure ${\rP}_0^{t,x}$ can also be viewed as a probability measure on the canonical space $C(\R^+,\R^d)$ endowed with the filtration $(\ccc_t)$ generated by the coordinate mappings. Let $\rP_0$ be the law on the canonical space of the standard Brownian motion. For $s\in [0,t)$, we have
\[
{\rm P}_0^{t,x}(d\rw)_{\mid \ccc_s}=\frac{p(t-s,x-\rw(s))}{p(t,x)} {\rm P}_0(d\rw)_{\mid \ccc_s}.
\]
We consider the energy of $\mu$ with respect to the process $(W_s)$ (see \cite{lg:hppt} for a precise description and definition). Thanks to \cite[Proposition 1.1]{lg:hppt} we have:
\[
\ce(\mu)
=2 \int_0^\infty  ds\;\rP_0 \left[\left(\iint_E \nu(dt,dx) p(t-s,x-\rw(s))/p(t,x)\right)^2\right]=2I(\nu).
\]
Now, using \cite[Proposition 5]{dlg:wtsbmbs}, we know there exists an additive  functional $A$ of the Brownian snake killed when its lifetime reaches $0$ such that:
\begin{itemize}
\item[(i)] For every  Borel function $F\geq 0$ on $\cw_{0,0}^*$,  
  $\N_{0,0}\left[\int_0^\infty  F(W_s)dA_s\right] =\int\mu(d\rw)F(\rw) $.
\item[(ii)] $\N_{0,0}[A_\infty ^2]=2\ce(\mu)$.
\end{itemize}
We deduce from (i) that the additive functional increases only when $\hat W_s\in \supp \nu\subset K$. Therefore, using the  Cauchy-Schwarz inequality, we get
\[
  \N_{0,0}[\cg^*\cap K\neq\emptyset] \geq \N_{0,0}[A_\infty >0] \geq \N_0[A_\infty ]^2/\N_0[A_\infty ^2].
\]
We get 
$\N_{0,0}[\cg^*\cap K\neq\emptyset] \geq [4I(\nu)]^{-1}$. Since the above inequality is true for any probability $\nu$ on $K$, we get that
\[
\N_{0,0}[\cg^*\cap K\neq\emptyset]\geq  4^{-1} \Capa (K)= 4^{-1} \capa (K).
\]

\noindent
\textbf{Proof} of Theorem \ref{th:ucap}.
Notice the application  defined on $\cb(E)$ by $T(A)=\N_{0,0}[\cg^*\cap
A\neq \emptyset]$ for $A\in \cb(E)$ is a Choquet capacity
(see~\cite[th\'eor\`eme 1]{c:fatc}). Since the capacity $\capa$ is  an inner capacity (see \cite[Theorem 12]{m:tcpf}), it is enough to prove the theorem for compact subsets of $E$. The result is then given by the  previous section (with $C_0=c_3^2 \gamma C$) and the above result.
\findemo

\section{Brownian range and support of $X_1$}

In this section, we first give an estimation for the hitting  probabilities of the support of $X_1$. Then we prove that the range of Brownian motion and the support of super-Brownian motion at fixed time are intersection equivalent. 

Let us fix $d\geq 2$. We denote by $\capa_{d-2}$ the usual Newtonian (logarithmic if $d=2$) capacity in $\R^d$:
\[
\capa_{d-2} (A)=\left[\inf \iint_{\R^d\times \R^d}\rho(dx)\rho(dy)\;h_{d-2}(\val{x-y})\right]^{-1},
\]
with $h_{\gamma}(r)=r^{-\gamma}$ if $\gamma>0$ and $h_0(r)=\log_+(1/r)$. The infimum is taken over all probability measures  $\rho$ on $\R^d$ such that $\rho(A)=1$. Let $B(0,h)$ be  the  open ball of ${\mathbb R}^d$ centered at $0$ with radius $h$.

\begin{proposition}
\label{prop:hitting}
Let $M>0$. There exist two positive constants $a$ and $b$ such that for any Borel set $A\subset B(0,1)$, for any finite measure $\mu$ on $B(0,1)$, with $(\mu,\ind)\leq M$, we have
\[
a (\mu,\ind) \capa_{d-2} (A)\leq \P_{0,\mu}[\supp X_1 \cap A \neq \emptyset]\leq b(\mu,\ind) \capa_{d-2} (A).
\]
\end{proposition}

\noindent
\textbf{Proof.} 
Let $A\subset B(0,2)$ be a Borel set.
Let $\nu$ be a probability measure on $E$ such that $\nu(\{1\}\times A)=1$. Then we have $\nu=\delta_{\{1\}}\times \rho$, where $\rho$ is a probability measure on $\R^d$ such that $\rho(A)=1$.
We get
\begin{multline*}
  I(\nu)=  \iint_{(0,1)\times \R^d}  dsdy\; p(s,y)\\
\iint_{A\times A} \rho(dx)\rho(dx')p(1-s,x-y)p(1-s,x'-y) {p(1,x)}^{-1}{p(1,x')}^{-1}.
\end{multline*}
Since $x,x'$ are in $B(0,2)$ and since $s\in(0,1)$ it is easy to see there exist two positive constants $a_1$ and $b_1$ (independent of $A$ and $\rho$) such that
\[
a_1 I(\nu)\leq \iint_{A\times A} \rho(dx)\rho(dx') \;h_{d-2}(\val{x-x'}) \leq b_1 I(\nu).
\]

\noindent
This implies that for any Borel set $A\subset B(0,2)$,
\[
a_1 \capa_{d-2} (A)\leq \capa (\{1\}\times A)\leq b_1 \capa_{d-2} (A).
\]
Since the capacity $\capa_{d-2}$ is invariant by translation, we get that for any Borel set $A\subset B(0,1)$, for any $x\in B(0,1)$,
\[
a_1 \capa_{d-2} (A)\leq \capa (\{1\}\times A_{x})\leq b_1 \capa_{d-2} (A),
\]
where $A_x=\left\{y; y-x\in A\right\}$.
We deduce from Theorem \ref{th:ucap} that 
\[
4^{-1} a_1 \capa_{d-2} (A)\leq \N_{0,x}[\cg^*\cap (\{1\}\times A)\neq \emptyset]  \leq C_0 b_1 \capa_{d-2} (A).
\]
Since $X_1=\sum_{i\in I}Y_1(W^i)$, where $\sum_{i\in I} \delta_{W^i}$ is a Poisson measure on
   $C(\R^+,\cw_{0})$ with intensity $\int\mu(dx)\N_{0,x}[\cdot]$, we have
\[
\P_{0,\mu}[\supp X_1 \cap A \neq \emptyset]=1-\expp{-\int \mu(dx)\N_{0,x}[\supp Y_1 \cap A \neq \emptyset]}.
\]
Notice that  $\N_{0,x}$-a.e., $\{1\}\times (\supp Y_1 \cap A )=\cg^*  \cap (\{1\}\times  A )$. Since $(\mu,\ind)<M$, we then easily get the result.
\findemo

Intersection-equivalence between random sets has been defined by Peres \cite{p:iebpbp}. Two random Borel sets $F_1$ and $F_2$ in ${\mathbb R}^d$ are intersection-equivalent in an open set $U$, if there exist positive constants $a$ and $b$ such that, for any Borel set $A\subset U$,
\[
a\, \rP [A\cap F_1]\leq \rP [A\cap F_2]\leq b\, \rP [A\cap F_1].
\]
If $\pi$ is a probability measure on $B(0,1)$, then we denote by $\rP_{\pi}$ the law of  a $d$-dimensional Brownian motion  $(B_t, t\geq  0)$ started with the law $\pi$. For $d\geq 3$ the range of Brownian motion is defined by $\range_B=\{B_t,t\geq 0\}$ in $\R^d$. For $d=2$, we also denote by $\range_B$ the set $\range_B=\{B_t,t\in [0,\xi]\}$, where $\xi$ is an exponential random variable of parameter $1$ independent of $(B_t,t\geq 0)$. 
						
\begin{corollary}
Let $M>0$. There exist two positive constants $a$ and $b$ such that for any Borel set $A\subset B(0,1)$, for any absolutely continuous probability measure $\pi$ on $B(0,1)$ with density bounded by $M$, for any finite measure $\mu$ on $B(0,1)$, with $(\mu,\ind)\leq M$, we have
\[
a (\mu,\ind) \rP_\pi[\range_B\cap A\neq\emptyset ]\leq \P_{0,\mu}[\supp X_1 \cap A \neq \emptyset]\leq b(\mu,\ind) \rP_\pi[\range_B\cap A\neq\emptyset ].
\]
\end{corollary}

\noindent
\textbf{Proof}.  This is a consequence of Proposition \ref{prop:hitting} and the fact that there exist two positive constants $a_2$ and $b_2$ such that for any Borel set $A\subset B(0,1)$, for any absolutely continuous probability measure $\pi$ on $B(0,1)$ with density bounded by $M$, 
\begin{equation*}
a_2 \capa_{d-2} (A)\leq \rP_\pi[\range_B\cap A\neq\emptyset ]\leq b_2 \capa_{d-2} (A)
\end{equation*}
(see for example \cite[Proposition 3.2]{p:iebpbp} for $d\geq 3$ and \cite{p:csf} for $d=2$).
\findemo

\section{Appendix}

\noindent
In this section, we give the proof of Lemma~\ref{lem:opL1}, which relies
on the properties of the Hermite polynomials. We first recall the definition and some properties
of those polynomials.

\subsection{Hermite polynomials}

\noindent
For  $n=(n_1,\cdots,n_d)\in\N^d$, we set $\val{n}=\sum_{i=1}^d n_i$, $n!=\prod_{i=1}^d n_i!$ and $\sum_{n\geq 0}=\sum_{i=1}^d\sum_{n_i=0}^\infty $. For $j\in \{1,\cdots,d\}$, let $\delta(j)$ be the element of $\N^d$ such that $\delta(j)_i=\delta_{i,j}$, the standard Kronecker symbol. If $z=(z_1,\cdots,z_d)$ is an element of $\R^d$, then we set $z^n=\prod_{i=1}^d z_i^{n_i}$. Let $(\cdot,\cdot)$ be the Euclidean product on $\R^d$.

The function $\varphi(z)=\expp{-[\val{z}^2-2(z,x)]/2}$ is an entire function defined on $\R^d$. We have
\begin{equation}
\label{eq:defH}
\expp{-[\val{z}^2-2(x,z)]/2}=\sum_{n\geq 0} \inv{n!} z^n He_n(x),
\end{equation}
where the $n$-th term $He_n(x)$  is a polynomial of $(x_1,\cdots,x_d)$ of degree $\val{n}$ called the $n$-th Hermite polynomial. Those polynomials can easily be expressed with the usual one dimensional Hermite polynomials $(He_k^{(1)}, k\in \N)$: 
$He_n(x)=\prod_{i=1}^d He_{n_i}^{(1)}(x_i)$, where $ x=(x_1,\cdots,x_d)$.

Now, let us recall some basic properties of the polynomials $He_n$. The following recurrence formula can be deduced from (\ref{eq:defH}) by derivating w.r.t. $z_i$: for all $n\in \N^d$ such that  $n_i>0$, 
\begin{equation}
   \label{eq:reccH}
He_{n}(x)=x_i He_{n-\delta(i)}(x)-(n_i-1)He_{n-2\delta(i)}(x),\quad \forall x\in \R^d,
\end{equation}
where by convention $He_{n-k\delta(i)}=0$ if $n_i-k<0$.
The derivative formula can be deduced from  \reff{eq:defH} by derivating w.r.t. $x_i$: for all $n\in \N^d$,
\begin{equation}
\label{eq:diffH}
\partial_i He_n=n_i He_{ n-\delta(i)}.
\end{equation}
We also recall the upper  bound for $He_n$ (see \cite[22.14.17]{as:hmff}): there exists a universal constant $1<c_0<2$ such that 
\begin{equation}
   \label{eq:boundH}
\val{He_n(x)}\leq c_0^d \sqrt{n!} \expp{\val{x}^2/4}\quad \text{for all}\quad x\in \R^d,\quad n\in \N^d.
\end{equation}
Using the definition of the Hermite polynomials, it is also easy to prove that:
\begin{equation}
   \label{eq:orthoH}
   \int dx\; p(t,x) He_n(x/\sqrt{t})He_m(x/\sqrt{t})=n!\prod_{i=1}^d \delta_{n_i,m_i}.
\end{equation}
It is also well known that the Hermite polynomials is a complete orthogonal system in $L^2(\R^d, \expp{-\val{x}^2/2} dx)$. Finally, standard arguments on Hilbert spaces show that if $f\in L^2(p)$ then 
\[
f(t,x)=\sum_{n\geq 0} f_n(t,x) =\sum_{n\geq 0} He_n(x/\sqrt{t}) g_n(t),
\]
where $g_n(t)=(n!)^{-1} \int dx\; p(t,x)He_n(x/\sqrt{t})f(t,x)$ and $g_n\in L^2((0,\infty ))$. Furthermore, we have
\begin{equation}
\label{eq:norme-de-f}
\norm{f}_{(p)}^2=\sum_{n\geq 0} n!\int_0^\infty dt\; g_n(t)^2.
\end{equation}
Since  $C_0^\infty ((0,\infty ))$ is dense in $L^2((0,\infty ))$, it is clear that the set $\ca$ of functions $f(t,x)= \sum_{n\geq 0} He_n(x/\sqrt{t}) g_n(t)$ where $g_n\in C_0^\infty ((0,\infty ))$ is non zero for a finite number of indices $n$, is dense in $L^2(p)$.
\label{com:Hermitecomplet}

\subsection{Proof of Lemma \ref{lem:opL1}}

\noindent
In a first step we prove there exist unique bounded extensions $\tilde\Lambda_1$, $\tilde\Lambda_{2,i}$ and  $\tilde\Lambda_{3,i}$ in $L^2(p)$ of the operators $\Lambda_1$, $\Lambda_{2,i}$ and  $\Lambda_{3,i}$ defined on $\ca$. Then in a second step we  check that the extensions $\tilde\Lambda_1$, $\tilde\Lambda_{2,i}$ and  $\tilde\Lambda_{3,i}$ and  the operators $\Lambda_1$, $\Lambda_{2,i}$ and  $\Lambda_{3,i}$, which are also defined on $C_0^\infty (E)$, agree on $C_0^\infty (E)$.

\textbf{First step}.
Let us compute $\Lambda (f)$ for very particular functions $f\in \ca$. Let $g\in C_0^\infty ((0,\infty ))$, $\alpha$ and $\beta$ be two positive reals such that $\supp g \in [\alpha,\beta]$, and $G(t)=\int_0^t ds\; g(s)$. For $n\in {\mathbb N}^d$, and $(t,x)\in E$, we set
\[
h_{n,g}(t,x)=He_n(x/\sqrt{t})t^{-\val{n}/2} g(t).
\]
Let us prove that
\begin{equation}
\label{eq:egLh}
\Lambda(h_{n,g})=h_{n,G}.
\end{equation}

For $z\in {\mathbb R}$, we introduce the function $\ch_{g,z}$ defined on $E$ by 
\[
\ch_{g,z}(t,x)=\sum_{n\geq 0} \inv{n!} z^n h_{n,g}(t,x) =g(t)\expp{-[\val{z}^2-2(x,z)]/2t}.
\]
Then we have
\begin{align*}
   \Lambda(\ch_{g,z})(t,x)
&=\inv{p(t,x)}\int_0^t ds \int dy\; p(t-s,x-y)p(s,y) \ch_{g,z}(s,y)\\
&=\inv{p(t,x)}\int_0^t ds \int dy\; p\left(\frac{s(t-s)}{t}, y-\frac{sx}{t}-\frac{(t-s)z}{t}\right)p(t,z-x) g(s)\\
&=\expp{-[\val{z}^2-2(x,z)]/2t}\int_0^t ds \; g(s)=\ch_{G,z}(t,x).\\
\end{align*}
Using \reff{eq:boundH},  Chapman-Kolmogorov equation, and that $\supp g \in[\alpha,\beta]$, we get
\begin{align*}
\Lambda(\val{h_{n,g}})(t,x)
&\leq \sqrt{n!} \;c_0^d p(t,x)^{-1}\iint_E ds dy\; p(t-s,x-y)p(s,y)\expp{\val{y}^2/4s} s^{-\val{n}/2} \val{g(s)}\\
&\leq \sqrt{n!} \;(\sqrt{2}c_0)^d p(t,x)^{-1}\int_0^t ds \;p(t+s,x)s^{-\val{n}/2} \val{g(s)}\\
&\leq \sqrt{n!} \;(\sqrt{2}c_0)^d \expp{\val{x}^2/4t}\int_0^t ds \;s^{-\val{n}/2} \val{g(s)}\\
&\leq \sqrt{n!} \;(\sqrt{2}c_0)^d \expp{\val{x}^2/4t}\norm{g}_\infty  (\beta-\alpha)\alpha^{-\val{n}/2}.
\end{align*}
The radius of the series $\sum a^k (k!\alpha^k)^{-1/2}$ is infinite. Thus for any $(t,x)\in E$ the series $\sum (n!)^{-1} z^n  \Lambda(\val{h_{n,g}})(t,x)$ are convergent. Fubini's theorem implies that
\[
\sum_{n\geq 0} \inv{n!} z^n  \Lambda(h_{n,g})(t,x)=\Lambda\left(\sum_{n\geq 0} \inv{n!} z^n h_{n,g}\right)(t,x)=\ch_{G,z}(t,x).
\]
Hence the two series $\sum (n!)^{-1} z^n  \Lambda(h_{n,g})(t,x)$ and $\sum (n!)^{-1} z^n  h_{n,G}(t,x)$ agree. Since their radius of convergence is positive (in fact infinite), we get that \reff{eq:egLh} is true.

Let us prove that $\Lambda_{2,i}$ has a bounded extension on $L^2(p)$.
We deduce from \reff{eq:diffH} that
\begin{align}
\nonumber
\Lambda_{2,i}(h_{n,g})(t,x)&=\frac{1}{2}\partial_{i,i}^{2}\Lambda(h_{n,g})(t,x)\\
\label{eq:lambda2i}
&=\inv{2t} n_i(n_i-1) He_{n-2\delta(i)}(x/\sqrt{t}) t^{-\val{n}/2}\int_0^t ds\;g(s).
\end{align}
Let us introduce $f\in \ca$, i.e. for $(t,x)\in E$, $f(t,x) =\sum_{n\geq 0} He_n(x/\sqrt{t})g_n(t)$, where $g_n\in C_0^\infty ((0,\infty ))$ and $g_n=0 $ except for a finite number of terms. By linearity, we have
\[
   \Lambda_{2,i}(f)(t,x)
= \sum_{n\geq 0} 2^{-1} n_i(n_i-1) He_{n-2\delta(i)}(x/\sqrt{t}) t^{-1-\val{n}/2}\int_0^t ds\; s^{\val{n}/2} g_n(s).
\]
Thus,  using \reff{eq:orthoH}, we have
\begin{align}
\nonumber 
  \norm{\Lambda_{2,i}(f)}_{(p)}^2
&= \sum_{n\geq 0} (n-2\delta(i))! 4^{-1}n_i^2(n_i-1)^2\int_0^\infty  dt \; t^{-2-\val{n}}\left[\int_0^t ds\; s^{\val{n}/2} g_{n}(s)\right]^2\\
\nonumber 
&\leq \sum_{n\geq 0} n! \frac{n_i(n_i-1)}{4}\frac{4}{(\val{n}+1)^2}\int_0^\infty  dt \;g_{n}(t)^2\\
\nonumber 
& \leq  \norm{f}_{(p)}^2,
\end{align}
\label{com:Hardy}
where we used the Hardy inequality: for $k>-1$,
\begin{equation*}
\int_0^\infty dt\; t^{-2-k}\left[\int_0^t s^{k/2} h(s)ds \right]^2
\leq \frac{4}{(k+1)^2}\int_0^\infty dt\;h(t)^2
\end{equation*}
for the first inequality and \reff{eq:norme-de-f} for the second one.
This means that $\Lambda_{2,i}$, defined on $\ca$, can be uniquely extended into a bounded operator $\tilde \Lambda_{2,i}$ from $L^2(p)$ into itself. The above inequality  implies $\normt{\tilde\Lambda_{2,i}}_{(p)} \leq 1$. 

For $i\in\{1,\cdots,d\}$, we set $\Lambda_{4,i}=\Lambda_{3,i}+2\Lambda_{2,i}$. Using \reff{eq:diffH} and  \reff{eq:reccH}, we deduce from \reff{eq:egLh} that 
\begin{multline}
\label{eq:lambda5i}
   \Lambda_{4,i}(h_{n,g})(t,x)\\
\begin{aligned}[b]
&=\left[{-(x_i/\sqrt{t})}\partial_i He_{n}(x/\sqrt{t})+n_i(n_i-1)He_{n-2\delta(i)}(x/\sqrt{t})\right]t^{-1-\val{n}/2}\int_0^tds\; g(s)\\
&=-n_i He_{n}(x/\sqrt{t})t^{-1-\val{n}/2}\int_0^tds\; g(s).
\end{aligned}
\end{multline}
Arguing as above, we get for $f\in \ca$,
\begin{align*}
\norm{\Lambda_{4,i}(f)}_{(p)}^2
&= \sum_{n\geq 0} n! n_i^2 \int_0^\infty  dt\; t^{-2-\val{n}}\left[\int_0^t ds\; s^{\val{n}/2}g_n(s)\right]^2\\
&\leq  \sum_{n\geq 0} n! n_i^2 \frac{4}{(\val{n}+1)^2} \int_0^\infty  dt\; g_n(t)^2\\
&\leq  4 \norm{f}_{(p)}^2.\\
\end{align*}
Thus the operators $\Lambda_{4,i}$ and  $\Lambda_{3,i}$, defined on $\ca$, can be uniquely extended in  bounded operators $\tilde\Lambda_{4,i}$ and  $\tilde\Lambda_{3,i}$ from $L^2(p)$ into itself. Furthermore we have $\normt{\tilde\Lambda_{4,i}}_{(p)}\leq 2$ and $\normt{\tilde\Lambda_{3,i}}_{(p)}\leq \normt{\tilde\Lambda_{4,i}}_{(p)}+2\normt{\tilde\Lambda_{2,i}}_{(p)}\leq 4 $. 

The proof concerning $\Lambda_1$  easily follows from the previous results. From \reff{eq:egLh}, we get 
\begin{multline*}
\Lambda_{1}(h_{n,g})(t,x)\\
=h_{n,g}(t,x)-\inv{2}\left[\val{n}He_n(x/\sqrt{t}) + \sum_{i=1}^d \frac{x_i}{\sqrt{t}} \partial_i He_n(x/\sqrt{t})\right]t^{-1-\val{n}/2}\int_0^t ds\; g(s).
\end{multline*}
Then using \reff{eq:lambda2i} and \reff{eq:lambda5i}, we get 
\[
\Lambda_1(h_{n,g})=\left[I+\inv{2}\sum_{i=1}^d \left[\Lambda_{4,i}+\Lambda_{3,i}\right]\right](h_{n,g})=\left[I+\sum_{i=1}^d \left[\Lambda_{4,i}-\Lambda_{2,i}\right]\right](h_{n,g}).
\]
This means that $\Lambda_1=I+\sum_{i=1}^d \left[\Lambda_{4,i}-\Lambda_{2,i}\right]$ on $\ca$. Hence $\Lambda_1$ can be uniquely extended in a bounded operator $\tilde\Lambda_1$ from $L^2(p)$ into itself and $\tilde\Lambda_1=I+\sum_{i=1}^d \left[\tilde\Lambda_{4,i}-\tilde\Lambda_{2,i}\right]$. We deduce that $\normt{\tilde\Lambda_1}_{(p)} \leq 1+3d$. 
\findemo

\textbf{Second step}. We first consider the operators $\Lambda_{3,i}$ for $i\in \{1,\cdots,d\}$.
To check that $\Lambda_{3,i}$ and $\tilde \Lambda_{3,i}$ agree on $C_0^\infty (E)$, it is enough to check that for $\varphi\in C_0^\infty (E)$, $\Lambda_{3,i}(\varphi)(t,x)=\tilde\Lambda_{3,i}(\varphi)(t,x)$ $dtdx$-a.e. Let $\varphi\in C_0^\infty (E)$. For $k\in \N$, we define,
\[
\varphi_k(t,x) =\sum_{\val{n}\leq k} He_n(x/\sqrt{t}) (n!)^{-1} \int_{\R^d} dy\; p(t,y)He_n(y/\sqrt{t})\varphi(t,y).
\]
The sequence $(\varphi_k,k\geq 0)$ converges in $L^2(p)$ to $\varphi$.

If $x\in \R^d$, $y\in \R$, $i \in \{1,\cdots,d\}$, we denote  by $z=\hat x^i_y$ the element of $\R^d$ such that $z_i=y$ and $z_j=x_j$ for $j\neq i$. Since $\Lambda_{3,i}(f)(t,x)=-t^{-1} x_i \partial_i \Lambda(f)(t,x)$ for $f\in \ca\cup C^\infty _0(E)$, we see that an integration by parts gives
\begin{equation}
   \label{eq:integbp}
\int_0^{x_i} dy\; \Lambda_{3,i}(f)(t,\hat x^i_y)= -t^{-1} x_i \Lambda(f)(t,x)+\int_0^{x_i} dy\; t^{-1} \Lambda (f)(t,\hat x^i_y).
\end{equation}
For short we write $P_i(f)$ for the operator $P_i(f)(t,x)=\int_0^{x_i} dy\; f(t,\hat x^i_y)$.
Let $R>0$ and  $T>\varepsilon>0$ be fixed. Let $Q=[\varepsilon,T]\times [-R,R]^d$. The heat kernel $p$ is bounded below and above on $Q$ by positive constant, say $c_Q$ and $C_Q$.
Using  Cauchy-Schwarz inequality we have
\[
\norm{\ind_Q P_i(f)}_{(p)}^2
\leq  C_Q R^2   \iint_Q dtdx \; f(t,x)^2 \leq  C_Q c_Q^{-1} R^2 \norm{f}_{(p)}^2.
\]
Thus the operator $\ind_Q P_i$ is continuous from $L^2(p)$ to $L^2(p)$. Thanks to Lemma \ref{lem:opL} and the above first step, we get that  the sequences $(\ind_Q P_i(\Lambda_{0}(\varphi_k)),k\geq 0)$ and $(\ind_Q P_i(\Lambda_{3,i}(\varphi_k)),k\geq 0)$ converge in $L^2(p)$ respectively to $\ind_Q P_i(\Lambda_{0}(\varphi))$ and $\ind_Q P_i(\tilde\Lambda_{3,i}(\varphi))$. Notice also that $(\ind_Q \Lambda (\varphi_{\sigma(k)}), k\geq 0)$ converges in $L^2(p)$ to $\ind_Q \Lambda (\varphi)$.
Thus,  there is a subsequence $(\sigma(k),k\geq 0)$ such that the sequences $(\ind_Q P_i(\Lambda_{0}(\varphi_{\sigma(k)})),k\geq 0)$, $(\ind_Q P_i(\Lambda_{3,i}(\varphi_{\sigma(k)})),k\geq 0)$ and $\left(\Lambda (\varphi_{\sigma(k)}),k\geq 0\right)$ converge  $dtdx$-a.e.  respectively to $\ind_Q P_i(\Lambda_{0}(\varphi))$, $\ind_Q P_i(\tilde\Lambda_{3,i}(\varphi))$ and $\Lambda (\varphi)$. Now 
\reff{eq:integbp} holds for $f=\varphi_{\sigma(k)}$, this means that for $(t,x)\in Q$,
\[
P_i (\Lambda_{3,i}(\varphi_{\sigma(k)}))(t,x)=-t^{-1} x_i \Lambda(\varphi_{\sigma(k)})(t,x)+P_i (\Lambda_{0}(\varphi_{\sigma(k)}))(t,x).
\]
Taking the limit we get that $dtdx$-a.e. in $Q$, 
\[
P_i (\tilde\Lambda_{3,i}(\varphi))(t,x)=-t^{-1} x_i \Lambda(\varphi)(t,x)+P_i (\Lambda_{0}(\varphi))(t,x).
\]
Since $R,T,\varepsilon$ are arbitrary, the above equality holds $dtdx$-a.e. in $E$.
Since \reff{eq:integbp} holds also for $f=\varphi$, we deduce that 
 $dtdx$-a.e., 
\[
\int_0^{x_i} dy \;  \Lambda_{3,i}(\varphi)(t, \hat x^i_y)=\int_0^{x_i} dy \; \tilde \Lambda_{3,i}(\varphi)(t, \hat x^i_y). 
\]
Hence we have $dtdx$-a.e., 
$\Lambda_{3,i}(\varphi)(t,x)=\tilde \Lambda_{3,i}(\varphi)(t,x)$.

The proofs concerning  the operators $\Lambda_1$ and $\Lambda_{2,i}$, for $i \in \{1,\cdots,d\}$, and their extensions follow the same ideas.
\findemo

\end{document}